\documentclass[12pt,reqno]{amsart}
\usepackage{amsthm,amsmath,amsfonts,amscd,amssymb,latexsym,exscale,xypic,verbatim}
\usepackage[mathscr]{eucal}
\input{defns.secants}

\begin{document}
\title{On Waring's problem for several algebraic forms}
\author{Enrico Carlini and  Jaydeep Chipalkatti}
\maketitle

\parbox{12cm}{\small
We reconsider the classical problem of representing a finite number of 
forms of degree $d$ in the polynomial ring over 
$n+1$ variables as scalar combinations of powers of linear forms. 
We define a geometric construct called a `grove', which, 
in a number of cases, allows us to determine the dimension of the space of 
forms which can be so represented for a fixed number of summands. 
We also present two new examples, where this dimension turns out to 
be less than what a na{\" \i}ve parameter count would predict. }
\medskip 
 
\section{Introduction} 
Waring's problem for algebraic forms is formulated in analogy with the 
number-theoretic version. Assume that $F_1,\dots,F_r$ are homogeneous forms of degree 
$d$ in variables $x_0,\dots,x_n$. We would like to find linear forms $Q_1,\dots,Q_s$, such 
that each $F_i$ is expressible as a linear combination of $Q_1^d,\dots,Q_s^d$. This 
problem, and especially the case $r=1$, has received a great deal of attention 
classically. Indeed, since the representation 
\begin{equation}
F = c_1 Q_1^d + \dots + c_s Q_s^d 
\label{formula.0} \end{equation}
is computationally easy to work with, geometric results about the hypersurface 
$F=0$ are sometimes more easily proved by reducing $F$ to such an  
expression by a linear change of variables. For instance, the classical texts of 
Salmon \cite{Salmon3, Salmon2} frequently use this device. 

Typically the forms $F_i$ were assumed general, and the goal of the enquiry was to find 
the smallest $s$ for which the problem is solvable. An elementary parameter count gives 
an expected value of $s$, which usually turns out to be correct. However, there 
are exceptional cases when the expected value does not suffice, and of course they
are the ones of more interest. Here we consider a more general version of 
the problem, i.e., we fix $s$ and ask for the dimension of the family of forms 
$(F_i)$ which can be so expressed. See \cite{Ger1, IK} for an overview of the 
problem. 

The formal set-up is as follows. 
Let $V$ be a $\complex$--vector space of dimension $n+1$, and consider the 
symmetric algebra $S = \bigoplus\limits_{d \ge 0} \Sym^d \, V$. Choosing a basis 
$\{x_0,\dots,x_n\}$ for $V$, an element in $S_d$ may be written as 
a degree $d$ form in the $x_i$. 

Fix two positive integers $r \le s$. 
Let $\uQ = \{ Q_1, \dots, Q_s \}$ denote a typical point of 
$\Sym^s (\P S_1)$, and consider the set 
\[ 
U_s = \{ \uQ : Q_1^d, \dots, Q_s^d \; \text{are linearly independent 
over $\complex$} \}.
\] 
This is an open set of $\Sym^s (\P S_1)$, and if $s \le \dim S_d$, then it is 
nonempty. (Indeed, if the $Q_i$ are chosen generally, then $Q_i^d$ are 
linearly independent--see \cite[p.~12 ff]{IK}.) Henceforth we assume 
$s \le \dim S_d$. 

Let $G(r,S_d)$ denote the Grassmannian of $r$-dimensional subspaces of $S_d$ and 
$\Lambda \in G(r,S_d)$ a typical point. 
Now consider the incidence correspondence 
$\Xi \subseteq G(r,S_d) \times U_s$, defined to be 
\begin{equation}
 \Xi  = \{ (\Lambda, \uQ): 
\Lambda \subseteq \vspan\, (Q_1^d,\dots,Q_s^d) \}.
\end{equation}
Let $\Sigma$ denote the image of the first projection 
$\pi_1:\Xi \lra G(r,S_d)$. The chief preoccupation of this 
paper is calculating the dimension of $\Sigma$. 

\begin{Remark} \rm 
In general $\Sigma$ may not be a quasiprojective variety. E.g., let $(n,d,r,s)= 
(1,3,1,2)$. A binary cubic $F$ lies in $\Sigma$, iff it is either a cube of a linear 
form, or has three distinct linear factors. Identify the set of cubes in $\P S_3$ 
with a twisted cubic curve $C$. Then its tangential developable 
$T_C$ (i.e.~the union of tangent lines to $C$) consists of forms 
which can be written as $Q_1^2Q_2, (Q_i \in S_1)$. Hence 
\[ \Sigma = (\P S_3 \setminus T_C) \cup C. \] 
In particular, the map $\pi_1|_\Xi$ may be dominant without being surjective. It is 
in general difficult to determine the smallest $s$ such that it is surjective, and we 
do not address this problem here. 
\end{Remark} 
\begin{Definition} \rm 
If $\uQ \in U_s$ and $\Lambda \subseteq \vspan(Q_1^d,\dots,Q_s^d)$, then 
$\uQ$ is called a polar $s$-hedron\footnote{If $n=2$,  we will 
of course say polar triangle, quadrilateral etc.} of $\Lambda$. 
\end{Definition}

Thus an element $\Lambda \in G(r,S_d)$ lies in $\Sigma$ iff it admits a
polar $s$-hedron. 
If $F_1,\dots,F_r$ span $\Lambda$, then 
we will speak of a polar $s$-hedron of the $F_i$. 

The projection $\pi_2: \Xi \lra U_s$ is a Grassmann bundle of 
relative dimension $r(s-r)$, hence 
$N_1 := \dim \, \Xi = sn + r(s-r)$. This is the number of parameters 
implicit in the right hand side of expression (\ref{formula.0}). Let 
$N_2 := \dim \, G(r,S_d) = r(\binom{n+d}{d}-r)$, then 
\begin{equation} 
\dim \Sigma \le \min \{N_1,N_2\}. 
\label{bound.dimsigma} \end{equation} 
We define the deficiency $\delta(\Sigma)$ as the difference 
$\min\{N_1,N_2\} - \dim \Sigma$. 
As we will see, positive deficiency is a rare phenomenon. 
A necessary condition for $\Sigma$ to be dense in $G(r,S_d)$ is 
$N_1 \ge N_2$, i.e., 
\begin{equation}
 s \ge \frac{r}{n+r} \binom{n+d}{d}. \label{bound.dense}
\end{equation}

If $\Sigma$ is dense in $G$, then the general fibre of $\pi_1: \Xi \lra 
\Sigma$ has dimension $N_1 - N_2$. An interesting case 
is $N_1 = N_2 = \dim \Sigma$, when a general $\Lambda$ admits finitely many 
polar $s$-hedra. But in very few cases we know how many. 

When $r=1$, a complete answer to the problem of calculating $\dim \Sigma$ is 
known. Using apolarity (or equivalently Macaulay--Matlis duality), the question is 
reduced to a calculation of the Hilbert function of general fat points in $\P^n$. 
The final theorem is due to Alexander and Hirschowitz \cite{AH}.
See \cite{Ger1,IK,RanestadSch} for further discussion and references. 

\begin{Theorem}[Alexander--Hirschowitz] 
Assume $r=1$ and $d \ge 3$. Then equality holds in (\ref{bound.dimsigma}) 
except when 
\[ 
(n,d,s) = (2,4,5),(3,4,9),(4,3,7) \; \text{or} \; (4,4,14). 
\] For all exceptions,  $\delta(\Sigma)=1$. 
\label{Theorem.AH} \end{Theorem} 

The case $r=1,d=2$ is anomalous, in the sense that $\Sigma$ is then almost 
always deficient. (See \cite[Ch.~22]{JoeH} for the exact calculation.)
Clebsch's discovery of the example $(2,4,5)$ (see \cite{Clebsch1}) was a surprise, 
as it showed that merely counting parameters was not sufficient to solve the problem. 
Thus a general planar quartic does not admit a polar pentagon, but a quartic which 
admits one (called a Clebsch quartic), admits at least $\infty^1$ of them. See 
\cite{DolgachevKanev} for some beautiful results on Clebsch quartics. 

In this paper we consider the case $r > 1$, which remains open in general. 
Terracini's paper \cite{Terra1} addresses this problem, but it is not easy to 
follow. 
We know of only four examples when $r > 1$ and 
(\ref{bound.dimsigma}) is not an equality, viz. 
\begin{equation} 
(n,d,r,s) = (2,3,2,5),(3,2,3,5),(3,2,5,6),(5,2,3,8),  
\label{exceptions} \end{equation}
with $\delta(\Sigma) =1$ in every case. The first two examples were classically known, 
see \cite{FranzLondon} for the first, and \cite[p.~353]{Darboux}, 
\cite{Frahm,Toeplitz} for the second. 
The last two were found by the authors using a computer search. 

The paper is organised as follows. In the next section we construct a 
morphism $\mu$ whose image is $\Sigma$. Then we differentiate the expression for 
$\mu$ to get a formula for the dimension of $\Sigma$ (see 
Theorem \ref{theorem.dimsigma}). This motivates the definition of a 
geometric construct called a `grove', which is, 
roughly speaking, a linear system of hypersurfaces with assigned 
singularities. In Theorem \ref{theorem.grove}, 
we reinterpret the codimension of $\Sigma$ as the dimension of a family of groves. 
In \S \ref{section.examples}, we give several examples to show how 
geometic arguments can used to calculate $\dim \Sigma$. 
In the last section, we try to prove the deficiency of the four examples above using 
this method. For the last example, we do not succeed entirely. 

\smallskip 
The first author would like to thank Prof.~K.~Ranestad for many fruitful conversations 
and Prof.~A.~V.~Geramita for suggesting the topic. 
The financial support of Oslo University, OMATS programme, and 
Pavia University was of great help during the preparation of this paper.

The second author would like to thank Prof.~A.~Iarrobino for his comments 
on a preliminary version of this paper. He gratefully acknowledges the financial 
support of Profs A.~V.~Geramita and L.~Roberts as well as Queen's University. 
We are also thankful to the authors of the computer algebra system Macaulay--2. 

\section{groves and the dimension of $\Sigma$}
\label{section.grove}
\subsection{An analytic representation of $\Sigma$}
Let $\Mat^\circ(1,r; S_d)$ be the set of matrices of size $1 \times r$ with 
entries in $S_d$, and columns independent over $\complex$. (Similar definitions 
are understood below.) Then 
$G(r,S_d)$ is the quotient $\Mat^\circ(1,r; S_d)/GL_{r}(\complex)$. If 
$C\Sigma$ denotes the inverse image of $\Sigma$ in $\Mat^\circ(1,r; S_d)$, then 
$\dim \Sigma = \dim C\Sigma -r^2$. 

Consider the morphism of varieties 
\begin{equation} \begin{aligned} 
\Mat^\circ(1,s;S_1) \times \Mat^\circ(s,r;\complex)& \stackrel{\mu}{\lra}
\Mat^\circ(1,r;S_d) \\ 
([Q_1,\dots, Q_s], A) & \lra 
[Q_1^d,\dots, Q_s^d]\, A = ([Z_1,\dots,Z_r]). 
\end{aligned} 
\label{morphism.matrices} \end{equation}
The image of $\mu$ is $C\Sigma$, hence $\dim C\Sigma$ is the rank of the 
Jacobian matrix of $\mu$ at a general point in the domain of 
$\mu$. 

We can use this setup for a machine computation of $\dim \Sigma$. Write 
\[ 
Q_i = \sum\limits_{j=0}^n q_{ij}x_j, \quad A=(\alpha_{ij}), 
\quad Z_k = \sum\limits_{|I|=d} z_{k,I}\, {\underline x}^I, \] 
where the $q,\alpha$ are indeterminates and $z_{k,I}$ polynomial functions in 
$q,\alpha$. The Jacobian 
\[ \partial (z_{k,I})}/{\partial (q,\alpha)
\]
is then easily written down, and in order to find its rank, we substitute
random numbers for the $q$ and $\alpha$.  
We programmed this in Macaulay-2 to search for deficient examples. 
The search shows that in the intervals below, 
there are no examples of deficiencies other than those already 
mentioned. 
\begin{itemize} 
\item 
$n =2$, $2 \le d \le 6$, all possible $r,s$ (recall that $s < \dim S_d$),
\item 
$n =3$, $2 \le d \le 3$, all possible $r,s$,
\item 
$n=4$, $d=2$, all possible $r,s$,
\item 
$n=4$, $d=3$, $r \le 14, s \le 23$,
\item 
$n=5$, $d=2$, all possible $r,s$,
\item 
$n=5$, $d=3$, $r \le 9, s \le 34$. 
\end{itemize} 

A.~Iarrobino pointed out that the deficient examples tend to occur for $s = n+2, n+3$, and when 
$N_1, N_2$ are close. However there are no further such examples in the following range:
\[ 2 \le n \le 10, 2 \le d \le 5, s = n+2, n+3, 
r = \lfloor \frac{ns}{\binom{n+d}{d} -s} \rfloor, 
\lceil \frac{ns}{\binom{n+d}{d} -s} \rceil.  \]

The source code for the Macaulay--2 routine is available upon request, for which 
the readers should contact the second author. 

\subsection{A formula for $\dim \Sigma$} 
We will now use the morphism $\mu$ to describe a formula for $\dim \Sigma$. 
Let $R = \bigoplus\limits_{d \ge 0} \Sym^d \, V^*$, so that $\P S_1 = \text{Proj}\, R$. 
If $X \subseteq \P S_1 (= \P^n)$ is a closed subscheme, then $I_X$ denotes 
its ideal and $I_X^{(2)}$ the second symbolic power of $I_X$. 

Let ${\underline Q} = \{Q_1,\dots, Q_s\}$ be a set of 
$s$ points in $\P S_1$. 
Given an $s \times r$ matrix $A = (\alpha_{ij})$ over $\complex$, we have a morphism 
\begin{equation} \begin{aligned}
\eta: \Mat(1,r; (I_{\uQ})_d) & 
\lra \bigoplus_{i=1}^s R_d/(I_{Q_i}^{(2)})_d  \\ 
[u_1,\dots,u_r] & \lra 
[\dots, \sum\limits_{j=1}^r \alpha_{ij}.u_j + (I_{Q_i}^{(2)})_d, \dots]_{1 \le i \le s}
\end{aligned} \label{formulaA} \end{equation} 

Of course $\eta$ depends on the choice of $A, \uQ$, but we 
will write $\eta_{A,\uQ}$ only if confusion is otherwise likely. 

\begin{Theorem}
With notation as above, assume that the points ${\underline Q}$ and the matrix 
$A$ are general. Then 
\begin{equation} \codim(\Sigma, G(r,S_d)) = \dim \ker(\eta).
\label{formula.dimsigma} \end{equation} 
\label{theorem.dimsigma} \end{Theorem}
The proof uses the classical notion of apolarity. We introduce the essentials, 
see e.g. \cite{EhRo, Ger1,GrYo,IK} for details. 

\subsection{Apolarity} 
Recall that 
\[ R=\bigoplus\limits_{d \ge 0} \Sym^d \, V^*, \quad  
 S=\bigoplus\limits_{d \ge 0} \Sym^d \, V. \] 
Let $\{ x_0,\dots, x_n \}$ and $\{\partial_0,\dots, \partial_n \}$ be 
the dual bases of $V$ and $V^*$ respectively. We interpret 
a polynomial $u(\partial_0,\dots,\partial_n)$ in $R$ as the differential 
operator 
$u(\frac{\partial}{\partial x_0}, \dots, \frac{\partial}{\partial x_n})$. 
Then we have maps $R_p \circ S_q \lra S_{q-p}$, and thus $S$ acquires the 
structure of an $R$--module. 

For a subspace $W \subseteq S_d$, let 
\[ W^\perp = \{ u \in R_d: u \circ F = 0 \; \; \text{for every $F \in W$} \},\]
which is a subspace of $R_d$, such that 
$\dim W^\perp + \dim W = \dim S_d$. In classical terminology, if 
$u \circ F =0$ and $\deg u \le \deg F$, then $u, F$ are said to be apolar to 
each other. Thus $W^\perp$ is the set of differential operators in $R_d$, which 
are apolar to all forms in $W$. 

In the following two instances $W^\perp$ can be concretely 
described (see \cite[Lemma 2.2]{IK}). Let $Q \in S_1$ be a nonzero 
linear form, or equivalently a point in $\P S_1$. 
\begin{enumerate}
\item[i.]{If $W = \text{span}\,(Q^d)$, then 
$W^\perp = (I_Q)_d$. }
\item[ii.]{If $W = \{ Q^{d-1}Q': Q' \in S_1 \}$, then 
$W^\perp = (I_Q^{(2)})_d$. \label{perp2}} 
\end{enumerate} 
\smallskip 

{\sc Proof of theorem} \ref{theorem.dimsigma}. 
 We will calculate the map on tangent spaces for the 
morphism $\mu$ in (\ref{morphism.matrices}). Fix a general point 
$([Q_1,\dots,Q_s],A)$. Given arbitrary forms $Q_1',\dots, Q_s' \in S_1$ and 
$B \in \Mat(s,r;\complex)$, we have 
\[ \begin{aligned} 
{} & \mu([Q_1+\epsilon Q_1',\dots,Q_s+\epsilon Q_s'], A + \epsilon B) - 
\mu([Q_1,\dots,Q_s],A) = \\
& \epsilon\{ [Q_1^d,\dots,Q_s^d]B + 
d[Q_1^{d-1}Q_1',\dots,Q_s^{d-1}Q_s']A \} + O(\epsilon^2). 
\end{aligned} \]
Hence the tangent space to $C\Sigma$ at the point $\mu([Q_1,\dots,Q_s],A)$ is described 
as 
\[ \begin{aligned} 
{\mathbb T} = \{[Q_1^d,\dots,Q_s^d]B & + [Q_1^{d-1}Q_1',\dots,Q_s^{d-1}Q_s']A: \\ 
& Q_1',\dots, Q_s' \in S_1, B \in \Mat(s,r;\complex) \}. 
\end{aligned} \] 
Now $\dim {\mathbb T} = \dim C\Sigma = \dim \Sigma + r^2$. 
Define maps 
\[ \begin{aligned} 
\alpha:  \Mat(1,s;S_1) & \lra \Mat(1,r;S_d)  \\ 
[Q_1',\dots,Q_s'] & \lra [Q_1^{d-1}Q_1',\dots,Q_s^{d-1}Q_s']A, \quad \text{and} \\ 
\beta:  \Mat(s,r;\complex) & \lra \Mat(1,r;S_d)  \\ 
B & \lra [Q_1^d,\dots,Q_s^d]B,  \end{aligned} \]
so that ${\mathbb T} = \image \alpha + \image \beta$. 
After dualising, we have a diagram 
\[ \diagram 
\Mat(1,r;R_d) \rto^{\alpha^*} \dto_{\beta^*} & \Mat(1,s;R_1) \\ 
\Mat(s,r;\complex) 
\enddiagram \] 

Now $\uu = [u_1,\dots,u_r] \in \ker \beta^* \iff$
for every $[F_1,\dots,F_r] \in \image \beta$, we have 
$u_i \circ F_i =0$ for all $i$. For any pair of indices $1 \le i_1,i_2 \le r$, 
one can certainly arrange $B$ such that $F_{i_1} = Q_{i_2}^d$. Thus 
$\uu \in \ker \beta^*$ iff each $u_i$ lies in 
$\cap_j \, \vspan(Q_j^d)^\perp = \cap_j \, (I_{Q_j})_d = (I_{\uQ})_d$. Hence 
$\ker \beta^* = \Mat(1,r; (I_{\underline Q})_d)$. 

By analogous reasoning, an element $\uu \in \ker \beta^*$ will be in 
$\ker \alpha^*$ iff it annihilates all elements in $\image \alpha$, i.e., 
iff for every $i$, the operator $\sum\limits_{j=1}^r \alpha_{ij}.u_j$ is apolar to 
$\{ Q_i.Q': Q' \in S_1\}$. Thus with the natural inclusion 
\[ \Mat(1,r;(I_{\uQ})_d) \subseteq \Mat(1,r;R_d), \] 
we have $\ker \eta = \ker \alpha^* \cap \ker \beta^*$. Finally 
\[ \begin{aligned} 
{} & \dim \ker \eta = \dim \ker \alpha^* + \dim \ker \beta^* 
- \dim (\ker \alpha^* + \ker \beta^*) \\ 
= & \, (r\,\dim R_d - \dim \image \alpha) + 
  (r\,\dim R_d - \dim \image \beta) \\ - & 
  (r \, \dim R_d - \dim (\image \alpha \cap \image \beta)) \\  
= & \, r\, \dim R_d - \dim (\image \alpha + \image \beta) \\ 
= & \, r\, \dim R_d - \dim {\mathbb T} 
= r\, \dim R_d - \dim \Sigma - r^2 
= \dim G(r,S_d) - \dim \Sigma. 
\end{aligned} \] 
The theorem is proved.  \qed 
\smallskip 

\noindent If $r=1$, then $\ker \eta = (I_{\uQ}^{(2)})_d$. 
Hence we recover the formula 
(see \cite[Theorem 6.1]{Ger1})
\begin{equation} \dim \Sigma = \dim\,  (R/I_{\underline Q}^{(2)})_d -1. 
\end{equation}

\begin{Remark} \rm Since $\dim  \ker \eta$ is upper semicontinuous in 
the variables $A, \uQ$ (see \cite[p.~125, exer.~5.8]{Ha})
\[ \dim \ker \eta \ge \codim \, \Sigma \ge \max\{0,N_2- N_1\} \] 
for \emph{any} choice of $A$ and ${\underline Q}$. Hence if the first and the 
last terms coincide for some choice, then it follows that 
$\Sigma$ is not deficient. 
\label{remark.uppersemi} \end{Remark} 
We will reformulate this theorem geometrically. In the sequel, assume that 
$\uQ = \{Q_1,\dots, Q_s\}$ are points in a \emph{fixed} copy of $\P^n$ 
($\P S_1$ if you will) and similarly $\up = \{p_1,\dots,p_s \}$ are points 
in $\P^{r-1}$. 

\begin{Definition} \rm 
A \emph{grove}\footnote{After some fitful experimentation, we decided to choose 
a name devoid of any mathematical associations.}
for the data $(\up, \uQ)$ consists of a triple $(\Gamma, L, \gamma)$ 
such that 
\begin{itemize}
\item{$\Gamma \subseteq \P H^0(\P^n,\O_{\P}(d))$ is a linear system of 
dimension (say) $t \le r-1$,}
\item{$L \subseteq \P^{r-1}$ is a linear space of dimension $r-(t+2)$ 
(thus defining a projection $\pi_L: \P^{r-1} -\ra \P^t$), and}
\item{$\gamma: \P^t \stackrel{\sim}{\lra} \Gamma$ is an isomorphism,} 
\end{itemize} 
satisfying the following conditions: 
\begin{itemize}
\item{all the $Q_i$ belong to the base locus of $\Gamma$,}
\item{for every $i$, either $p_i \in L$ or the hypersurface 
$\gamma \circ \pi_L(p_i)$ is singular at $Q_i$.}
\end{itemize} \label{defn.grove} \end{Definition}
We denote the collection of all groves by $\grove \, (\up, \uQ)$. 
\begin{Remark} \rm 
To make the definition of $\pi_L$ canonical, identify 
$\P^t$ with the set of linear subspaces of dimension $r-(t+1)$ containing $L$, and 
then let $\pi_L(p) = \overline{L p}$. If $t=r-1$, then $L$ is taken as empty and 
$\pi_L$ the identity map. (In the applications, almost always this will be 
the case.) If $L = \emptyset$, then $\Gamma$ is an $(r-1)$-dimensional system of 
degree $d$ hypersurfaces passing through $\uQ$, such that $\gamma(p_i)$ is 
singular at $Q_i$. 

If $r=1$, then necessarily $t=0, L = \emptyset$ and all $p_i$ are the 
same point. Then a grove is a solitary hypersurface of degree $d$ singular 
at all $Q_i$. 
\label{remark.grove} \end{Remark} 

For the next proposition, we identify $\P^{r-1}$ with $\P \, \Mat(1,r; \complex)$. 
If $A \in \Mat(s,r; \complex)$ is a matrix with no zero rows, then we identify its 
$i$-th row as the point $p_i \in \P^{r-1}$. 

\begin{Proposition} Fix points $Q_1,\dots,Q_s$ in $\P^n$. Then with identifications 
as above, we have a bijection $\P(\ker \eta_{A,\uQ}) \simeq \grove \, (\up, \uQ)$. 
\label{proposition.grove} \end{Proposition}
\demo 
Let $\uu = [u_1,\dots,u_r]$ be a nonzero element of $\ker \eta$. Let 
$\Gamma$ be the linear system generated by the $u_i$, and 
\[ 
L = \{ X \in \Mat(1,r;\complex): [F_1,\dots,F_r]\, X^t = 0 \}. 
\] 
Then $\pi_L$ appears as the map 
\[\P \, \Mat(1,r;\complex) \, -\ra \P (\Mat(1,r;\complex)/L) (= \P^t). \] 
Define 
\[ \gamma: \P^t \stackrel{\sim}{\lra} \Gamma, \quad  
   X + (L) \lra [u_1,\dots,u_r]\, X^t. \] 
By hypothesis, the form $\alpha_{i1}u_1 + \dots \alpha_{ir}u_r$ lies in 
$(I_{Q_i}^{(2)})_d$. Hence, unless it is identically zero 
(i.e., $p_i = [\alpha_{i1},\dots, \alpha_{ir}] \in L$), the 
hypersurface it defines (which is $\gamma \circ \pi_L(p_i)$) is singular at $Q_i$. 

Alternately, given a grove $(\Gamma, L, \gamma)$, assume that
$\Gamma$ is defined by $W \subseteq H^0(\P^n, \O_\P(d))$. Then $\gamma$ 
induces an isomorphism $\hat \gamma: \Mat(1,r)/L \lra W$ (well-defined upto 
a global scalar). 
Now if $u_i = \hat \gamma([0,\dots,1,0, \dots])$ (the $1$ in $i-$th place), 
then $\uu \in \ker \eta$. This defines the bijection. 
\qed 

The next result follows directly from Theorem \ref{theorem.dimsigma}. 
Nearly all subsequent results are based on this reformulation. 
\begin{Theorem} Let points $p_1, \dots, p_s \in \P^{r-1}$ and 
$Q_1, \dots, Q_s \in \P^n$ be chosen generally. Then $\Sigma$ has codimension 
$c$ in $G(r,S_d)$ if and only if, there are exactly 
$\infty^{c-1}$ groves for $(\up, \uQ)$. In particular, 
$\Sigma$ is dense in $G(r, S_d)$ if and only if, 
the points $(\up, \uQ)$ do not admit a grove. 
\label{theorem.grove} \end{Theorem}

In the paper of Terracini cited above, he states something which resembles 
the last statement in the theorem. 
Unfortunately, neither his statement nor the argument leading to it 
are clear. 

In the case $r=1$, we recover the criterion of Ehrenborg and Rota \cite[Theorem 4.2]{EhRo}. 
\begin{Corollary}[Ehrenborg, Rota] A general form in $S_d$ cannot be written as a 
sum of $d$-th powers of $s$ linear forms if and only if, given general points 
$Q_1,\dots,Q_s$ in $\P^n$, there exists a hypersurface of degree $d$ singular at 
all of them. 
\end{Corollary} 
Consider the collection 
\[ \grove^\circ \, (\up, \uQ) 
 = \{(\Gamma, L, \gamma): \text{$L$ contains none of the $p_i$} \}. 
\]
\begin{Lemma} 
Assume that the points $(\up, \uQ)$ are general. Then 
$\grove^\circ \, (\up, \uQ)$ 
is a nonempty Zariski open subset of $\grove \, (\up,\uQ)$.  
\label{lemma.generalgrove} \end{Lemma} 
Hence for purposes of calculating 
$\dim \Sigma$, we can assume that our groves lie in $\grove^\circ$. 

\demo 
Let $\grove_i \in \P(\ker \eta)$ be the open set of groves where $p_i \not \in L$, 
then $\grove^\circ = \cap_i \, \grove_i$. Thus $\grove^\circ$ fails to be dense only if 
some $\grove_i$ is empty. But then by symmetry (here is where the generality is used) 
each $\grove_i$ is empty, implying that every $L$ contains all the $p_i$. Since 
the set $\up$ spans $\P^{r-1}$ (recall $s \ge r$), this is impossible. 
\qed 

From Remark \ref{remark.uppersemi}, we know that 
\[ \dim \, \grove \, (\up, \uQ) \ge \codim \, \Sigma -1 \ge \max \{0, N_2 - N_1 \}-1, \] 
for any choice of points $(\up, \uQ)$. If the end terms are equal for some configuration 
of points, then $\Sigma$ is not deficient.

\section{examples} 
\label{section.examples} 
In this section we give a rather large number of examples illustrating the 
use of Theorem \ref{theorem.grove}. All the results follow the same plan: we choose 
specific values of $(n,d,r,s)$, then calculate the dimension of $\grove$ and  hence that of 
$\Sigma$. The choice of quadruples $(n,d,r,s)$ does not follow any definite 
pattern, but we have given examples which we think are geometrically interesting. 
Some of the results proved here are known, and the novelty lies 
in the method used to obtain them. 

We refer to \cite{JoeH} for the miscellaneous geometric facts needed. 
We mention two which will be used 
frequently. Recall that a set of points in $\P^n$ is said to be in 
linearly general position if any subset of $m$ points ($m \le n+1$) 
is not contained in a $\P^{m-2}$. 
\begin{itemize} 
\item 
Given two sequences 
$\{A_1, \dots, A_{n+2}\}, \{B_1, \dots, B_{n+2} \} \subseteq \P^n$ in linearly general 
position, there is a unique automorphism $\gamma$ of $\P^n$, such that 
$\gamma(A_i) = B_i$ for 
all $i$. 
\item 
Given $n+3$ points of $\P^n$ in linearly general position, 
there is a unique rational normal curve passing through all of them. 
\end{itemize} 

For every case treated in this section, $\dim \Sigma$ will coincide with the expected value 
$\min\{N_1, N_2\}$. The deficient examples are the subject of the next section. 

The following result should be classically known, but we have 
been unable to trace a reference. 

\begin{Theorem}
If $n=1$, then $\Sigma$ is not deficient for any $d,r,s$. 
\end{Theorem}
\demo Let $Q_1,\dots,Q_s$ and $A = (\alpha_{i,j})$ be as above. 
Consider the composite map of vector 
bundles on $\P^1$: 
\[ \rho_A: \{ \O_{\P^1}(d H - \sum Q_i)\}^{\oplus r}  \lra 
\{ \O_{\P^1}(d H)\}^{\oplus r}  
\stackrel{\tilde\eta}{\lra} \bigoplus\limits_{i=1}^s \O_{2Q_i}(dH) \] 

Here $H$ denotes the hyperplane divisor on $\P^1$. The map on the 
left is the canonical 
inclusion, and the one on the right is induced by $A$. On local sections, 
\[ \begin{aligned} 
({\tilde \eta}[u_1,\dots,u_r])_i = 
\sum\limits_{j=1}^r \alpha_{ij}u_j, \;\; & \text{modulo functions vanishing to} \\ 
& \text{order at least $2$ at $Q_i$.} 
\end{aligned} \] 

The map $H^0(\P^1,\rho_A)$ is identical to $\eta$ in formula (\ref{formulaA}). 
Hence if $\E = \ker \rho_A$, then $h^0(\E) = \codim \, \Sigma$. 
The image of $\rho_A$ is the skyscraper sheaf 
\[ \bigoplus\limits_i \ker \, (\O_{2Q_i}(d H)  \lra 
\O_{Q_i}(d H)) = 
\bigoplus\limits_i \O_{Q_i}(d H - Q_i) \] with degree $s$, hence 
$\E$ is a rank $r$-vector bundle of degree $\epsilon = r(d-s)-s$.

Now specialise $A$ to the following matrix: 
write $s = r \alpha  + \beta$, with $0 \le \beta \le r-1$ and let 
\[ A^t = \left[ B_1 | \dots |B_{r-\beta}\,| C_{r+1-\beta}| \dots | C_r \right], 
\; \text{where} \] 
\begin{itemize} 
\item 
the $B_i$ (resp. $C_i$) are blocks of size $r \times \alpha$ 
(resp. $r \times (\alpha+1)$), 
\item 
each $B_i$ or $C_i$ is made of all $1$'s in the $i$-th row and zeros elsewhere. 
\end{itemize} 
Then $\E$ splits as a direct sum 
\begin{equation} \O_{\P^1}(d-\alpha-s)^{\oplus (r-\beta)} \oplus  
\O_{\P^1}(d-\alpha -s -1)^{\oplus \beta}. \label{twistsofE} \end{equation}
Now $N_1 = s + r (s-r)$ and $N_2 = r (d-r + 1)$, so $N_2 - N_1 
= \epsilon + r$. If $N_2 \le N_1$, then all twists in (\ref{twistsofE}) are  
negative, so $h^0(\E) = 0$. If $N_2 > N_1$, then all twists are at least $-1$, 
so $h^0(\E) = N_2 - N_1$. In either case 
$\codim \, \Sigma = \max\{0,N_2-N_1\}$, hence by Remark \ref{remark.uppersemi}
we are through. 
\qed 

\begin{Remark} \rm 
Fix points $\uQ$, and think of $\E$ as moving in a family parametrised by 
$A$. By Grothendieck's theorem,  $\E$ splits into a direct sum of line bundles. 
The point of the theorem is that if $A$ is general, then its splitting 
type is balanced, i.e., it deviates from the 
sequence $(\deg \E/\rank \E, \dots,\deg \E/\rank \E)$ as little 
as possible. Once the splitting type is known, $h^0(\E)$ is known. 
\end{Remark} 

\begin{Example} \rm This example might give some insight into the construction 
of $A$. Let $r=3, s=7$, so $\alpha = 2, \beta =1$. Then 
\[ A^t = \left[ \begin{array}{ccccccc} 
1 & 1 & 0 & 0 & 0 & 0 & 0 \\ 
0 & 0 & 1 & 1 & 0 & 0 & 0 \\ 
0 & 0 & 0 & 0 & 1 & 1 & 1 \end{array} \right] 
\] 
and $\tilde{\eta}([u_1,u_2,u_3]) = [u_1,u_1,u_2,u_2,u_3,u_3,u_3]$. 
Thus a local section $[u_1,u_2,u_3]$ will lie in $\ker \rho_A$ 
iff $u_1$ (resp.~$u_2$ and $u_3$) vanishes doubly at $Q_1,Q_2$ 
(resp.~ at $Q_3,Q_4$ and $Q_5,Q_6,Q_7$). Hence 
$\E$ is a direct sum 
\[ \begin{aligned} 
{} & \O_{\P^1}(dH - Q_1-Q_2 - \sum Q_i) \oplus 
\O_{\P^1}(dH - Q_3 - Q_4 - \sum Q_i) \oplus \\ 
& \O_{\P^1}(dH - Q_5 - Q_6 - Q_7 - \sum Q_i). 
\end{aligned} \] 
\end{Example} 

Henceforth we use the same notation for a form $F \in S_d$ and 
the hypersurface in $\P R_1$ which it defines. 

\begin{Proposition} 
Two general plane conics have a unique polar 
triangle. $(N_1 = N_2 =8.)$
\label{prop.2223} \end{Proposition} 
Firstly we will show that $\dim \, \Sigma(2,2,2,3) = 8$. Choose 
general points $p_1,p_2,p_3 \in \P^1, Q_1,Q_1,Q_3 \in \P^2$, and let 
$(\Gamma, L, \gamma) \in \grove^\circ$ be a grove. 
Since there is no conic singular at all $Q_i$, 
$\dim \Gamma =1$ and $L = \emptyset$. 
Now $\gamma(p_1)$ must be the line pair $Q_1Q_2 + Q_1Q_3$ and similarly for 
other $p_i$. Since any two elements $\gamma(p_i), \gamma(p_j)$ span 
$\Gamma$, all the three lines $Q_iQ_j$ are in the base locus 
of $\Gamma$. This is absurd, hence there is no such grove. 

Consequently, two general conics $F_1, F_2$ admit at least one polar 
triangle--say $\{Q_1,Q_2,Q_3 \}$.\footnote{These $Q_i$ are 
unrelated to those in the previous paragraph. By the nature of our 
deductions, the $Q_i$ lead a 
double life: they are alternately linear forms and points.}
Now the pencil generated by the $F_i$ contains a member belonging to 
$\vspan(Q_1^2,Q_2^2)$, and this member must be singular at the 
point $Q_1 \cap Q_2$. Hence the points $Q_i \cap Q_j$ must be the vertices 
of the three line pairs contained in the pencil. 
This gives a geometric construction of the 
polar triangle and simultaneously shows that it is unique: 

Let $F_1, F_2$ intersect in $\{Z_1, \dots, Z_4 \}$. 
Let $A_1$ be the point of intersection of the lines $Z_1Z_2, Z_3Z_4$, and 
similarly $A_2 = Z_1Z_3 \cap Z_2 Z_4, A_3 = Z_1Z_4 \cap Z_2Z_3$. Define  
lines $Q_1 = A_2A_3, Q_2 = A_1A_3, Q_3 = A_1A_2$. Then 
$\{Q_1,Q_2,Q_3 \}$ is the required triangle. \qed 

\begin{Proposition} 
Four general plane conics $F_1,\dots,F_4$ have a unique polar 
quadrilateral. $(N_1 = N_2 =8.)$
\label{prop.2244} \end{Proposition} 
\demo Firstly let us show that $\Sigma\, (2,2,4,4)$ is dense in $G(4,S_2)$. 
Let $p_1, \dots, p_4 \in \P^3, Q_1, \dots, Q_4 \in \P^2$ be chosen 
generally, and $(\Gamma, L, \gamma) \in \grove^\circ(\up, \uQ)$. Since there is no conic 
singular at all $Q_i$, we must have $\dim \Gamma =1$. Then 
$\Gamma$ is the pencil of conics through $\uQ$, which has no members singular at any 
$Q_i$. This precludes any possibility of defining $\gamma$. 

Thus four general conics $F_1, \dots, F_4$ admit at least one polar quadrilateral, 
say $\{Q_1, \dots, Q_4 \}$. We may assume that $Q_i$ are in linearly  
general position. 
Let $A = [\alpha_0, \alpha_1, \alpha_2]$ be the point of intersection of 
the lines $Q_1,Q_2$. (Thus as an element of $\P R_1$, 
$A = \alpha_0 \partial_0 + \alpha_1 \partial_1 + \alpha_2 \partial_2$ upto a 
scalar). By hypothesis, 
\[ F_1 = c_1 Q_1^2 + \dots + c_4 Q_4^2, \quad \text{for some constants $c_i$}. 
\] 
Operate by $A$ on the equality above, then 
\[ A \circ F_1 = \sum 2c_i Q_i(A) Q_i. \] 
Now $Q_1(A) = Q_2(A) = 0$, hence 
$A \circ F_1$ (the polar line of $F_1$ with respect to $A$) belongs to the pencil 
generated by lines $Q_3, Q_4$. An identical argument applies to 
all $F_i$, hence we deduce that the four lines 
$A \circ F_1, \dots, A \circ F_4$ are concurrent at the point $Q_3 \cap Q_4$. 
The line $A \circ F_i$ has equation 
\[ 
\frac{\partial F_i}{\partial x_0}(A) \, x_0 + \frac{\partial F_i}{\partial x_1}(A) \, x_1 + 
\frac{\partial F_i}{\partial x_2}(A) \, x_2 = 0, 
\] 
hence the Jacobian matrix 
$J = \partial(F_1,\dots,F_4)/\partial(x_0,x_1,x_2)$, 
has rank at most two at $A$. 

Now consider the locus $X = \{ \rank J \le 2 \} \subseteq \P^2$. 
It is easily seen that 
$X$ must be a finite set. Hence we have a Hilbert-Burch (or Eagon-Northcott) 
resolution 
\[ 0 \lra S(-4)^3 \lra S(-3)^4 \lra S \lra S/I_X \lra 0. \] 
From the resolution (or the Porteous formula), 
we have $\deg X=6$. By the argument above $X$ contains the points $Q_i \cap Q_j$, 
so it can contain no others. 

We claim that this forces the polar quadrilateral to be unique. Indeed let 
$M_1$ be a side of such a quadrilateral. The argument shows that 
$M_1$ must contain three of the points from $X$. This is impossible 
unless $M_1$ coincides with one of the $Q_i$. 
\qed 

\begin{Proposition} The variety $\Sigma\, (2,2,3,3)$ has 
dimension $6$. $(N_1 = 6, N_2 = 9.)$
\label{prop.2233} \end{Proposition} 

\demo Let $p_1,p_2,p_3, Q_1,Q_2, Q_3$ be general points in $\P^2$. 
We will show that $\up, \uQ$ admit exactly $\infty^2$ groves. Let 
$(\Gamma,L,\gamma) \in \grove^\circ$. Let $G_1$ be the line pair 
$Q_1Q_2 + Q_1Q_3$, 
and similarly for $G_2, G_3$. Evidently each $G_i$ belongs to $\Gamma$, 
hence $ \Gamma = \vspan(G_1, G_2, G_3)$ and $L = \emptyset$. 
Thus the only moving part of the grove is $\gamma$, and 
$\grove^\circ$ is isomorphic to the variety 
\[ \{ \gamma: \P^2 \stackrel{\sim}{\lra} \Gamma \; \text{such that 
$\gamma(p_i) = G_i$ for $i=1,2,3$} \}. \] 
Fix a point $Z \in \P^2$ such that $p_1, p_2,p_3,Z$ are in linearly general position. 
Then $\gamma$ is entirely determined by 
$\gamma(Z)$, so $\grove^\circ$ is isomorphic to an open set of 
$\P^2$. \qed 

We will frequently use B{\'e}zout's theorem in the following form: if a 
hypersurface of degree $d$ intersects a curve of degree $e$ in a scheme of length 
$> de$, then it must contain the curve. In such a circumstance we will loosely 
say that the hypersurface contains at least $de +1$ points of the curve. 

\begin{Theorem}[Sylvester's pentahedral theorem]
A general cubic surface in $\P^3$ has a polar 
pentahedron. $(N_1 = N_2 = 19.)$
\end{Theorem} 
The statement says that $\Sigma \, (3,3,1,5)$ is dense in $\P^{19}$, and 
it is covered by the Alexander--Hirschowitz theorem. We give a short geometric proof. 

\demo Choose general points $Q_1,\dots,Q_5$ in $\P^3$ and assume that a cubic 
$F$ is singular at all of them. Choose a sixth general point $Z$ and let $C$ be the 
unique twisted cubic through $Q_1,\dots,Q_5,Z$. Since $F$ contains at least $10$ 
points of $C$ (counting each $Q_i$ as two points), it must 
contain $C$ by B{\'e}zout's theorem. This implies the absurdity that $F$ contains 
a general point of $\P^3$. Hence there is no such $F$ and the claim is proved. 
\qed 

In \cite{Sylvester1}, Sylvester asserted that a general quaternary 
cubic has a unique polar pentahedron, and adduced some cryptic remarks in 
support. See \cite{RanestadSch} for a proof of the uniqueness. 

The next result is a direct generalisation of Proposition \ref{prop.2223}. 
\begin{Theorem} The variety $\Sigma(n,2,2,n+1)$ is dense in $G(2,S_2)$, 
moreover two general quadrics in $\P^n$ admit a unique 
polar $(n+1)$-hedron. $(N_1 = N_2 = n^2 + 3n-2.)$
\end{Theorem}
\demo Choose general points $p_1,\dots,p_{n+1} \in \P^1$, and 
$Q_1,\dots,Q_{n+1} \in \P^n$ and let $(\Gamma,L,\gamma) \in \grove^\circ$. 
There is no quadric singular at all $Q_i$ (since the singular locus of 
a quadric is a linear space, and the $\uQ$ are not contained in any proper 
linear subspace), hence $\dim \Gamma=1$ and $L = \emptyset$. The quadric 
$\gamma(p_i)$ contains at least three points of the line $Q_iQ_j$ 
(viz.~$Q_i$ twice and $Q_j$), so it must contain the line. Since any two 
quadrics $\gamma(p_i),\gamma(p_j)$ span $\Gamma$, it follows that all the 
lines $Q_iQ_j$ lie in the base locus of $\Gamma$. 

Let $F \in \Gamma$ and $F(-,-)$ its associated bilinear form. 
By what we have said, $F(Q_i+\lambda Q_j,Q_i+\lambda Q_j) = 0$ for all 
$\lambda \in \complex$, hence $F(Q_i,Q_j) =0$. Since the $Q_i$ span $\P^n$, we 
have $F \equiv 0$. This is absurd, so $(\up, \uQ)$ do not admit a grove. 

The proof of uniqueness is similar to Proposition \ref{prop.2244}. Let 
$F_1, F_2$ be general quadrics in $\P^n$ admitting a polar 
$(n+1)$-hedron $\{Q_1, \dots, Q_{n+1} \}$. Define points 
$A_i = \bigcap\limits_{j \neq i} Q_j \in \P^n$ for $1 \le i \le n+1$. For any 
$i$, the polar hyperplanes $A_i \circ F_1, A_i \circ F_2$ coincide, hence 
the Jacobian matrix
$J = \partial(F_1,F_2)/\partial(x_0,\dots,x_n)$ must have  rank one at each 
$A_i$. Now let $X = \{ \rank J \le 1 \}$, and use Hilbert-Burch together with Porteous 
to show that $X = \{A_1, \dots, A_{n+1}\}$. Then $Q_i$ is uniquely determined 
as the linear span of the points $A_j$ ($j \neq i$). \qed 
\begin{Remark} \rm 
Before proceeding we record a small construction for later use. 
Let $C$ be a twisted cubic in $\P^3$, and let 
$\Psi \subseteq \P H^0(\P^3, \O_\P(2))$ be the two-dimensional linear 
system of quadrics containing $C$. For every $x \in C$, there is a unique 
quadric (say $\psi_x$) in $\Psi$ singular at $x$. Thus we have an imbedding 
\[ \tau: C \lra  \Psi, \quad x \lra \psi_x.  \] 
Its image $\tau(C)$ is a smooth conic in $\Psi$. 

This notation will come in force only when we explicitly refer to this 
remark. Otherwise $C,\Psi$ etc may have unrelated meanings. 
\label{remark.twistedcubic} \end{Remark} 

The following technical result will be useful later. 
\begin{Lemma} Let $f,v: \P^1 \lra \P^2$ be two morphisms. 
Assume that $f$ is birational onto its image which is a curve of degree $m$, 
and $v$ is an imbedding onto a smooth conic. Assume moreover, that 
there are $m+2$ 
points $\lambda_1, \dots, \lambda_{m+2}$ in $\P^1$, such that 
$f(\lambda_i) = v(\lambda_i)$ for all $i$. Then $v=f$. 
\label{lemma.peculiar} \end{Lemma} 

\demo Choose a coordinate $x$ on $\P^1$ such that $\lambda_{m+2} = \infty$. 
We may choose coordinates on $\P^2$ such that $v(x) = [1,x,x^2]$. Then 
$f(x) = [A_0, A_1, A_2]$, such that $A_i$ are polynomials in $x$ with no 
common factor and $\deg A_i \le m$. By hypothesis, $f(\infty) = [0,0,1]$, hence 
$\deg A_2 > \deg A_1, \deg A_0$. In particular, $\deg A_0 \le m-1$. Now the 
polynomial $A_1 - xA_0$ (which is of degree $\le m$), vanishes for $m+1$ 
values $\lambda_1, \dots, \lambda_{m+1}$, hence it vanishes identically. But then 
$\deg A_0 \le m-2$. By the same argument, $A_2 - x^2A_0$ vanishes identically, 
hence $[A_0, A_1, A_2] = [1,x,x^2]$. \qed 

\begin{Remark} \rm 
If $C$ is a curve isomorphic to $\P^1$ and $A_1,\dots,A_4$ distinct 
points on $C$, then $\langle A_1, A_2, A_3, A_4 \rangle_C$ will denote 
their cross-ratio as calculated on $C$. Of course, it depends on the choice $C$, 
for instance four points in $\P^2$ have different cross-ratios as calculated on different 
smooth conics passing through them. 
\end{Remark}

In 1870, Darboux claimed that the case $\Sigma(3,2,4,6)$ is 
deficient (see \cite[p.~357]{Darboux}). 
In \cite{Terra1}, Terracini states (without proof) that Darboux's claim is 
wrong, and in fact there is no  deficiency. Here we substantiate Terracini's 
statement. 

\begin{Proposition} The variety $\Sigma(3,2,4,6)$ is dense in $G(4,S_2)$.  
$(N_1 = 26, N_2 = 24.)$ \end{Proposition}  

\demo Choose general points $(\up, \uQ)$ as usual, where $\up$ and 
$\uQ$ lie in nominally distinct copies of $\P^3$. We can identify the copies 
in such a way that the following holds: 
$p_1,\dots,p_6, Q_1,\dots, Q_6$ are in the same
$\P^3$ so that $p_i = Q_i$ for $1 \le i \le 5$ and 
$p_6, Q_6$ are distinct general points. 

Let $(\Gamma, L, \gamma) \in \grove^\circ(\up, \uQ)$, and let $C$ be the unique
twisted cubic  through the $\uQ$. The quadric $\gamma \circ \pi_L(p_i)$ 
intersects $C$ in at least seven points, so must contain $C$. 
Hence necessarily $\gamma \circ \pi_L(p_i) = \psi_{Q_i}$ in the notation of 
Remark \ref{remark.twistedcubic}. Thus $\Gamma = \Psi$ and $L$ is a point in $\P^3$. 
Let $\P^2_{\langle L\rangle}$ be the set of lines through $L$ 
(cf.~Remark \ref{remark.grove}), so we have a map 
$\P^2_{\langle L\rangle} \stackrel{\gamma}{\lra} \Psi$ . 

Now there are two maps $C \lra \P^2_{\langle L\rangle}$, namely 
$\pi_L$ and $\gamma^{-1} \circ \tau$. The image of the latter (say $D$) is a 
smooth conic. Moreover, $\deg \image(\pi_L) \le 3$ and the two maps coincide 
on points $p_1,\dots,p_5 \, (=Q_1, \dots, Q_5)$. Hence by Lemma \ref{lemma.peculiar}, 
they must be the same. 
In particular, $\deg \pi_L(C) =2$ which is only possible if $L$ is a point \emph{on} $C$. 
We claim that $\pi_L(p_6) = \pi_L(Q_6)$. Indeed, since $\pi_L$ is 
an isomorphism on $C$, 
\[ \begin{aligned} 
{} & \langle \pi_L(p_1), \pi_L(p_2),\pi_L(p_3), \pi_L(p_6) \rangle_D \\
=  \; & \langle \psi_{Q_1}, \psi_{Q_2},\psi_{Q_3},\psi_{Q_6} \rangle_{\tau(C)}, \\
=  \; & \langle Q_1, Q_2,Q_3,Q_6 \rangle_C \\
=  \; & \langle \pi_L(Q_1), \pi_L(Q_2),\pi_L(Q_3), \pi_L(Q_6) \rangle_D 
\end{aligned} \]
which shows the claim. This implies that the chord $LQ_6$ (in case $L \neq Q_6$) 
or the tangent to $C$ at $L$ (in case $L = Q_6$) passes through $p_6$. 
Now for a fixed $Q_6$, the chords $\{LQ_6\}_{L \in C}$ fill 
only a surface in $\P^3$. Hence if we choose $p_6$ off this surface, then 
no such configuration can exist. 
Thus general points $(\up, \uQ)$ do not admit a grove, which proves the proposition. 
It follows that four general space quadrics have $\infty^2$ polar $6$-hedrons. 
\qed 

\begin{Proposition} 
The variety $\Sigma \, (4,2,2,4)$ has dimension $20$. 
$(N_1 = 20, N_2 = 26.)$
\label{prop.4224} \end{Proposition} 
\demo 
Choose general points $p_1,p_2,p_3,p_4 \in \P^1$ and 
$Q_1,\dots,Q_4 \in \P^4$. We will show that there are exactly 
$\infty^5$ groves for these data. Let $\Pi$ denote the $3$--space spanned 
by the $Q_i$, and choose $(\Gamma, L, \gamma) \in \grove^\circ$. 
If $\dim \Gamma =0$, then $\Gamma$ is $\Pi$ doubled, and $L$ any point on $\P^1$. 
Since this is only a one-dimensional family, we may assume 
$\dim \Gamma =1, L = \emptyset$. 

Each of the quadrics $\gamma(p_i), \gamma(p_j)$ contains three points of the line 
$Q_iQ_j$, hence contains the line. Since these quadrics span $\Gamma$, 
all six lines $Q_iQ_j$ are in the base locus of $\Gamma$. This forces 
$\Pi$ to be in the base locus. Hence there exists a unique $2$-plane 
$\Psi_\Gamma \subseteq \P^4$, such that 
\[ \Gamma = \Pi \, \text{(fixed component)} + \; 
\text{pencil of $3$-planes through $\Psi_\Gamma$}. \] 

This leads to the following construction: let $\Psi \in G(3,5)$ be a $2$--plane 
in $\P^4$ away from the $Q_i$ and let $\psi_1,\dots, \psi_4$ be the 
$3$--planes through $\Psi$ containing the points $Q_1, \dots, Q_4$ respectively. 
Now we have a rational map 
\[ f: G(3,5)\,  - \ra \P^1, \quad 
\Psi \lra \langle \psi_1,\psi_2,\psi_3,\psi_4 \rangle. \] 
It is easy to see that $f$ is nonconstant, hence dominant. Now if $\Psi$ 
belongs to the fibre $f^{-1}(\langle p_1,p_2,p_3,p_4 \rangle)$, then (and only then) 
we can define 
\[ \gamma: \P^1 \stackrel{\sim}{\lra} \Gamma, 
\quad p_i \lra \Pi + \overline{\Psi Q_i} \;\; \text{for $i=1,\dots,4$.} \] 
Thus $\grove^\circ$ is birational to the fibre 
$f^{-1}(\langle p_1,p_2,p_3,p_4 \rangle)$, which is five dimensional. 
\qed 

\begin{Proposition}[London \cite{FranzLondon}]
The variety $\Sigma (2,3,3,6)$ is dense in $G(3, S_3)$, i.e., 
three general plane cubics admit a polar hexagon. $(N_1 = N_2 = 21.)$
\label{prop.2336} \end{Proposition}
\noindent London's proof is laborious, and it may be doubted whether 
it meets modern standards of rigour. 

\demo 
It is enough to show that for \emph{some} configuration $(\up, \uQ)$, there 
is no grove (cf.~Remark \ref{remark.uppersemi}). 

Let $p_1,\dots,p_6$ be general points in $\P^2$. Fix a line $M$ in 
$\P^2$, take $Q_4,Q_5, Q_6$ to be general points on $M$ and 
$Q_1, Q_2, Q_3$ general points in $\P^2$ (away from $M$). Let 
$(\Delta, L, \delta)$ be\footnote{The change in notation is of course deliberate.}
in $\grove^\circ(\up, \uQ)$. Since there is no 
cubic singular at all $Q_i$, $\dim \Delta \ge 1$. Now $L$ is either a 
point or empty, in either case the cubics $\delta \circ \pi_L(p_i)$ ($i=4,5,6$) 
must span $\Delta$. Now any of them intersects $M$ in at least four points, 
so must contain it. Thus $M$ lies in the base locus of $\Delta$, and 
$\Delta = M$ (fixed component) $+ \,\Gamma$, where $\Gamma$ is a system of conics through 
$Q_1,Q_2,Q_3$. Since each of $Q_1,Q_2,Q_3$ is a singular point of some member of $\Gamma$, 
we have 
\[ \Gamma = \vspan(G_1,G_2, G_3), \] 
following the notation used in the proof of Proposition \ref{prop.2233}. 
In particular $L = \emptyset$. Composing the isomorphism 
$\Delta \lra \Gamma$ with $\delta$, we have an isomorphism 
$\gamma: \P^2 \lra \Gamma$ such that 
$(\Gamma, \emptyset, \gamma)$ is a grove of conics 
for $(p_1,p_2,p_3,Q_1,Q_2,Q_3)$. 
Think of $\gamma$ as belonging to the two--dimensional family 
in Proposition \ref{prop.2233}. 

For $i=4,5,6$, if $\lambda_i \subseteq \Gamma$ be the line consisting of 
conics passing through $Q_i$, then by hypothesis 
$\gamma(p_i) \in \lambda_i$. But the conditions 
$\gamma(p_4) \in \lambda_4, \gamma(p_5) \in \lambda_5$ determine $\gamma$ 
uniquely. (To see this point, choose coordinates on 
$\P^2, \Gamma$ such that 
\[ p_1,G_1 = [1,0,0], \; p_2,G_2 = [0,1,0], \; p_3,G_3 = [0,0,1], \; p_4=[1,1,1] \] 
and $\lambda_4$ has line coordinates $[1,1,1]$. Then the matrix of $\gamma$ is 
diagonal, say equal to 
$\left[ \begin{array}{ccc} 
a & 0 & 0 \\ 
0 & b & 0 \\ 
0 & 0 & c \end{array} \right]$. 
Since $\gamma(p_4) \in \lambda_4$, we have $a+b+c=0$, and 
$\gamma(p_5) \in \lambda_5$ forces another independent condition. 
But then the matrix is uniquely determined upto a scalar.)

We conclude that the grove $(\Delta,L,\delta)$ is entirely 
determined by the data $p_1,\dots,p_5, Q_1,\dots,Q_5$. This is absurd, since 
one can certainly choose $p_6,Q_6$ such that $\gamma(p_6) \notin \lambda_6$. 
Hence $(\up, \uQ)$ do not admit a grove. 
\qed 

After a lengthy analysis, London concludes that three general cubics admit 
\emph{two} polar hexagons. It would be worthwhile to re-examine his argument. 
We hope to take it up elsewhere. 
\section{exceptional cases}
In this section we will construct groves 
showing that $\Sigma$ is deficient for the four quadruples mentioned 
in the introduction. Part I of our construction for the case $(3,2,3,5)$ is 
built on a hint in Terracini \cite{Terra1}. The rest we believe to be new. 
As we confessed earlier, we have only partial success in the last case. 
\begin{Theorem} 
The variety $\Sigma(2,3,2,5)$ has codimension $1$ in $G(2,S_3)$. 
$(N_1 = N_2 = 16.)$
\end{Theorem} 
\noindent Part I (construction of the grove). 
Choose general points $0,1,\infty,\alpha,\beta$ in $\P^1$, and 
$Q_0,Q_1,Q_\infty,Q_\alpha,Q_\beta$ in $\P^2$. Let $C$ be the unique smooth conic 
through the $\uQ$. The proposed construction is as follows: let 
$Z$ be a point in $\P^2$ and 
\[ \Gamma = C \, \text{(fixed component)} + 
\text{pencil of lines through $Z$.} 
\] 
Then we define 
\[ 
\gamma: \P^1 \stackrel{\sim}{\lra} \Gamma, \quad 
 \star \lra C + \text{ line $ZQ_\star$} \; \; \text{for \, $\star = 0,1,\infty,\alpha,\beta$.} 
\]
Of course, for such a $\gamma$ to exist, the cross--ratios must agree. 
Hence the position of $Z$ is crucial. 

Let $D_\alpha$ denote the unique smooth conic through $Q_0,Q_1,Q_\infty,Q_\alpha$ 
such that $\langle Q_0,Q_1,Q_\infty,Q_\alpha \rangle_{D_{\alpha}} = \alpha$. 
Similarly, let $D_\beta$ be the unique smooth conic 
through $Q_0,Q_1,Q_\infty,Q_\beta$ such that 
$\langle Q_0,Q_1,Q_\infty,Q_\beta \rangle_{D_{\beta}} = \beta$.  

Let $D_\alpha \cap D_\beta = \{ Q_0,Q_1,Q_\infty,Z\}$. Since $Z$ lies on 
$D_\alpha$, we have 
$\langle ZQ_0, ZQ_1, ZQ_\infty, ZQ_\alpha \rangle = \alpha$ and similarly 
for $\beta$. Hence the sequences 
\[ \{0,1,\infty,\alpha,\beta \} , \quad 
\{ ZQ_0, ZQ_1, ZQ_{\infty}, ZQ_\alpha, ZQ_\beta \}, \] are projectively equivalent. 
This ensures that $\gamma$ is well-defined and we are through. 

Part II (uniqueness of the grove). 
In part I, we have shown that $\dim \grove \ge 0$ for general $(\up,\uQ)$, hence 
this is true of any $(\up,\uQ)$. If we show that the grove is unique 
for some configuration, it will follow that $\dim \grove =0$ for general 
$(\up,\uQ)$. 

Let $M,N$ be distinct lines in $\P^2$.  
Choose general points $Q_0,Q_1,Q_\infty$ on $M$ and $Q_\alpha,Q_\beta$ on $N$. Let  
$0,1,\infty,\alpha,\beta$ be general points of $\P^1$, and assume that  
$(\Gamma, L, \gamma)$ is a grove for these data. Since there is 
no cubic singular at all $Q_i$,  
$\dim \Gamma =1, L = \emptyset$. By B{\'e}zout, the cubics $\gamma(0),\gamma(1)$
contain $M$, hence $M$ is in the base locus of  $\Gamma$. Now $\Gamma \setminus M$ 
is a pencil of conics, which, by the same argument on 
$\gamma(\alpha), \gamma(\beta)$, contains $N$ in its base locus. 
Hence 
\[ \Gamma = M + N \; \text{(fixed components)} 
+ \; \text{pencil of lines through a point $Z$}. \] 
Now  map $\P^1$ to $M$, by sending $0,1,\infty$ to $Q_0,Q_1,Q_\infty$ respectively, and  
via this map, think of $\alpha,\beta$ as points on $M$. Then $Z$ is forced to 
be the point of intersection of the lines $\alpha.Q_\alpha, \beta.Q_\beta$. 
The grove is thus uniquely determined. The theorem is proved. \qed 

\begin{Remark} \rm There is a simple explanation for the deficiency of 
$\Sigma$. Let $F_1,F_2$ be two plane cubics admitting a polar 
pentagon $\{Q_1, \dots, Q_5\}$. 
Since $\vspan(F_1, F_2) \subseteq \vspan(Q_1^3,\dots,Q_5^3)$, 
we deduce that the six partial derivatives $\partial F_i / \partial x_j$ 
($i=1,2, j=0,1,2$) lie in $\vspan(Q_1^2,\dots,Q_5^2)$. Hence they must be linearly 
dependent, which amounts to a nontrivial algebraic condition on the $F_i$. 
It is easy to write this condition as the vanishing of a $6 \times 6$ determinant 
whose entries are functions in the coefficients of $F_i$ (see \cite{FranzLondon}). 
\end{Remark} 

For the next two theorems the notation of Remark \ref{remark.twistedcubic} will 
remain in force.
\begin{Theorem} 
The variety $\Sigma(3,2,3,5)$ has codimension $1$ in $G(3,S_2)$. 
$(N_1 = N_2 = 21.)$ \label{theorem.3235} \end{Theorem} 
\demo 
Choose general points $p_1,\dots,p_5$ in $\P^2$ and 
$Q_1,\dots,Q_5$ in $\P^3$. Let $E$ be the smooth conic through the $p_i$, 
and consider the imbedding 
\[ E \lra \Sym^3 E , \quad p \lra 3p.\]
Abstractly  $\Sym^3 E \simeq \P^3$, hence there is a unique 
isomorphism $\beta: \Sym^3\, E \lra \P^3$, such that 
$\beta(3p_i) = Q_i$. Let $C$ be the twisted cubic obtained as the image of 
the composite $ E \lra \Sym^3\, E \stackrel{\beta}{\lra} \P^3$.  

Part I (construction of the grove). 
Let $\Gamma = \Psi$ (in the notation of Remark \ref{remark.twistedcubic}) and define 
\[ \gamma: \P^2 \stackrel{\sim}{\lra} \Gamma, \quad 
p_i \lra \psi_{Q_i}\;\; \text{for $1 \le i \le 4$.}
\] 
The sequences $\{p_1,\dots,p_5\} \subseteq E, 
\{\psi_{Q_1},\dots,\psi_{Q_5}\} \subseteq \tau(C)$ are such 
that the cross-ratios of any two corresponding subsequences of four points 
are equal. Hence $\gamma(E) = \tau(C)$ and $\gamma(p_5) = \psi_{Q_5}$, 
implying that $(\Gamma, \emptyset, \gamma)$ is a grove. 

\smallskip 

Part II (uniqueness of the grove).  We now show that $\grove^\circ = \grove^\circ(\up,\uQ)$
is a singleton set. The plan of the proof 
is to choose a general element $g \in (\Gamma, L, \gamma) \in \grove^\circ$, and 
then to show that the generality forces it to be the same as the grove 
constructed above. By construction, the functions 
\[  \grove^\circ \lra \dim L, \quad \grove^\circ \lra \rank \gamma \circ \pi_L(p_i) 
= \rho_i \] 
are respectively upper and lower semicontinuous. 
(We mean the rank of $\gamma(-)$ as a quadric in $\P^3$.) 
Let $U_i \subseteq \grove^\circ$ be the open set where $\rho_i$ is maximal, and 
let $g \in \cap \, U_i$. By symmetry, all $\rho_i$ equal the same number $\rho$, 
which is either $2$ or $3$. (It cannot be $1$ since no plane can contain all 
$Q_i$.)

Case $\rho =3$. Each quadric $S_i = \gamma \circ \pi_L(p_i)$ is a cone with 
its vertex at $Q_i$. Then 
\[ S_i \cap S_j = \text{(line $Q_iQ_j$)} \cup C_{ij}, \] 
where $C_{ij}$ is a twisted cubic through $Q_1,\dots,Q_5$. 
For any three indices $i,j,k$, the 
quadrics $S_i,S_j,S_k$ span $\Gamma$. Hence the base locus of $\Gamma$ equals 
$ S_i \cap S_j \cap S_k$,  which is set-theoretically just 
$C_{ij} \cap C_{ik} \cap C_{jk}$. 

Assume that the base locus of $\Gamma$ is zero dimensional, then it is supported only 
on ${Q_1,\dots,Q_5}$ (since two twisted cubics can have at most 
five points in common). Moreover the $S_i$ intersect transversally at each $Q_j$, so 
each $Q_j$ is a reduced point of the base locus. This is a contradiction, since 
by B{\'e}zout, the base locus is a scheme of length $8$. Hence the base locus is positive 
dimensional, i.e., all $C_{ij}$ are the same twisted cubic $C$. 

It follows that $\Gamma = \Psi$ in the notation of Remark 
\ref{remark.twistedcubic}. Then $\gamma(p_i)$ must equal $\psi_{Q_i}$ for 
each $i$, which determines $\gamma$ uniquely. 
Hence $\grove = \grove^\circ$ is a singleton set whose ``general'' element 
is the one we have constructed in Part I. 
\smallskip 

Case $\rho =2$. We will show that this case is impossible. Each 
$S_i = \gamma \circ \pi_L(p_i)$ consists of two planes both of which 
pass through $Q_i$. 
We claim that the base locus of $\Gamma$ contains a line. Indeed 
$S_1,S_2$ contain the line ${Q_1Q_2}$. If it is not in the base locus, then none of 
the other $S_i$ can contain it. Then $S_3$ is the union of 
planes $Q_1Q_3Q_4 \cup Q_2Q_3Q_5$, and similarly for $S_4,S_5$. But then 
$S_3,S_4,S_5$ contain the line $Q_1Q_3$  (and $Q_2Q_5$), so it is in 
the base locus. 

Let $U_{ij} \subseteq \grove^\circ$ be the open set of groves which do not contain 
the line $Q_iQ_j$ in their base locus. If (say) $U_{12}$ is nonempty, then by symmetry 
each $U_{ij}$ is nonempty. Then a general element 
$g \in \cap \, U_{ij}$ (which by hypothesis has $\rho =2$) 
can contain none of the lines, which is a contradiction. Thus $U_{ij} = \emptyset$, 
implying that a general $\Gamma$ must contain all ten lines $Q_iQ_j$ in the base locus. 
This is surely impossible, hence $\rho \neq 2$. The proof of the theorem is 
complete. 
\qed 

\begin{Example} \rm 
Now let $\Pi$ be a plane in $\P^3$, and $Q_1, \dots, Q_4$ general points in $\Pi$. 
Choose $Q_5 \in \P^3$ generally (away from $\Pi$) and $p_1, \dots, p_5$ general 
points in $\P^2$. We know that this configuration admits a grove, let $(\Gamma, L, \gamma)$ 
be one. The quadric $\gamma \circ \pi_L(p_1)$ is 
singular at $Q_1$, moreover by B{\'e}zout, it contains the four lines 
$Q_1Q_i$. This would be impossible if the quadric were of rank $3$, hence it 
must contain $\Pi$. The same argument applies to $Q_2,Q_3,Q_4$, hence 
$\Gamma = \Pi$ (as fixed component) $+$ a system of planes through $Q_5$. 
But then no member of $\Gamma$ can be singular at $Q_5$, hence $\pi_L(p_5)$ 
is undefined, i.e., $L = p_5$. The base locus of the system of planes is a line, 
say $N$. This leads to the following construction: let 
$\P^2_{\langle Q_5 \rangle}$ denote the variety of lines through $Q_5$, and 
define 
\[ f: \P^2_{\langle Q_5 \rangle} - \ra \P^1, \quad 
N \lra \langle NQ_1, NQ_2, NQ_3, NQ_4 \rangle.  \] 
Let $\lambda$ denote the cross-ratio 
$\langle p_5p_1, p_5p_2, p_5p_3, p_5p_4 \rangle$. Now if $N \in f^{-1}(\lambda)$, 
then (and only then) we can define a grove as above. Thus 
$\grove \, (\up, \uQ)$ is a one-dimensional family, which demonstrates the 
upper-semicontinuity of $\dim \grove$. 
Secondly, Lemma \ref{lemma.generalgrove} fails for this set of points. 
\end{Example} 

\begin{Remark} \rm 
The following explanation of the deficiency is given by 
Salmon (\cite[vol.~I, Ch.~IX, \S 235]{Salmon3}). Let $F_1, F_2, F_3$ be quadratic forms 
in $x_0,\dots,x_3$. Introduce indeterminates $a,b,c$, and let $G = a F_1 + b F_2 + c F_3$. 
Then the discriminant $\Delta$ of $G$ (as a quadratic form in the $x_i$) is a 
quartic in $a,b,c$. Now by choosing $F_i$ generally, $\Delta$ can be made equal to 
any planar quartic. However, if we assume that the $F_i$ admit a polar 
pentahedron, then $\Delta$ is necessarily a L{\"u}roth quartic 
(see \cite{DolgachevKanev}). Since L{\"u}roth quartics form a hypersurface in 
$\P S_4$, this imposes an algebraic condition on $F_i$. 
\end{Remark}

\begin{Theorem} The variety $\Sigma(3,2,5,6)$ has 
codimension $3$ in $G(5,S_2)$. $(N_1 = 23, N_2 = 25.)$
\end{Theorem} 
\demo 
Choose general points $p_1,\dots,p_6 \in \P^4$ and 
$Q_1,\dots,Q_6 \in \P^3$. Let $C$ be the unique twisted cubic through 
the $Q_i$. There is a unique imbedding 
\[ \alpha: C \lra \P^4, \quad \alpha(Q_i) = p_i \; \; 
\text{for $1 \le i \le 6$.}  
\] 
Part I (construction of the groves). 
We will show that there are at least $\infty^2$ groves for 
these data. Let $\Gamma = \Psi$ in the notation of Remark 
\ref{remark.twistedcubic}. 
Let $L$ be a chord or a tangent of the rational normal quartic 
$\alpha(C)$. Let $\P^2_{\langle L \rangle}$ denote the collection of $2$--planes 
in $\P^4$ containing $L$, and 
\[ 
\pi_L : \P^4 -  \ra \P^2_{\langle L \rangle}, \quad p \lra \overline{L p} 
\] the natural projection. Now $\pi_L$ is defined everywhere on $\alpha(C)$, and  
$\pi_L(\alpha(C)) = D_L$ is a smooth conic in $\P^2_{\langle L \rangle}$. The 
sequences 
$\{Q_1,\dots,Q_6\} \subseteq C, \{ \pi_L(p_1), \dots, \pi_L(p_6) \} \subseteq D_L$ 
are such that any corresponding subsequences of four points have the 
same cross-ratio. Define 
\[ \gamma_L: \P^2_{\langle L \rangle} \stackrel{\sim}{\lra} \Gamma, 
\quad \pi_L(p_i) \lra \psi_{Q_i} \; \; 
\text{for $1 \le i \le 4$.} \] 
By what we have said, $\gamma_L(D_L) = \tau(C)$ and 
$\gamma_L \circ \pi_L(p_i) = \psi_{Q_i}$ for $i=5,6$. Thus 
$(\Gamma, L, \gamma_L)$ is a two-dimensional family of groves. 
\smallskip 

Part II (bounding the dimension of $\grove$). 
We will show that we have already constructed a dense set of possible groves. 
Let $(\Gamma, L, \gamma) \in \grove^\circ(\up, \uQ)$. 
Each $\gamma \circ \pi_L(p_i)$ contains at least 
seven points of $C$, hence contains $C$ by B{\'e}zout. Thus $C$ is in the base locus 
of $\Gamma$, i.e., $\Gamma \subseteq \Psi$.  Since $\Psi$ contains a unique element 
singular at $p_i$, $\Gamma = \Psi$  which in turn implies $\dim L = 1$. 
Let $\P^2_{\langle L \rangle}$ have the same meaning as above, so 
we have an isomorphism $\P^2_{\langle L \rangle} \stackrel{\gamma}{\lra} 
\Psi$. 

Now there are two maps $\alpha(C) \lra \P^2_{\langle L\rangle}$, namely 
$\pi_L$ and $\gamma^{-1} \circ \tau \circ \alpha^{-1}$. The image of the latter is a 
smooth conic. Moreover, $\deg \image(\pi_L) \le 4$ and the two maps coincide 
on points $p_1,\dots,p_6$. Hence by Lemma \ref{lemma.peculiar}, 
they must be the same. 
In particular, $\deg \image(\pi_L) = 2$ which is only possible if $L$ intersects 
$\alpha(C)$ twice. This implies that the grove belongs to the family 
constructed above. The theorem is proved. \qed 

The case $(5,2,3,8)$ is perhaps more surprising than the rest of the exceptions. By 
counting parameters, we expect three general quadrics in $\P^5$ to have 
$\infty^1$ polar octahedrons, but they do not have any. 

\subsection{The Segre-Gale transform}
Consider the variety $(\P^1)^8$ with the group $\text{Aut}(\P^1)$ acting 
componentwise. Let 
$U \subseteq (\P^1)^8$ be the open set of semistable points and 
$Y = U /\text{Aut}(\P^1)$ the GIT quotient. 

In the sequel, $\sigma: \P^1 \times \P^2 \lra \P^5$ denotes the Segre imbedding. 
Let $\underline{A} = A_1, \dots, A_8 \in \P^1, \up = p_1, \dots, p_8 \in \P^2$ be 
general points, and $C$ the unique rational normal quintic through the eight points 
$\sigma(A_i \times p_i)$. Choosing  an isomorphism $\alpha: C \lra \P^1$, we get 
a point 
\[ \underline{B} = 
(\alpha \circ \sigma(A_1 \times p_1), \dots, \alpha \circ \sigma(A_8 \times p_8))
\in Y, \] 
which we call the Segre-Gale transform of $(\underline{A}, \up)$. 
The passage via $\alpha$ between eight general points in $\P^5$ and 
eight points in $\P^1$ is an instance of the Gale transform--see 
\cite{DolgaOrtland,EisenPopescu}. 

\begin{Lemma} Fix eight general points $\up \in \P^2$. Then the rational map 
\[ 
\omega(\up): Y - \ra Y , \quad {\underline A} \lra {\underline B} 
\] 
is dominant. (The reader should check that it is well-defined.) 
\end{Lemma} 

\demo This is a direct computation using coordinates (and was done in Maple). Let 
\[ 
{\underline A} = (0,1,\infty,a_1,\dots,a_5), \quad 
p_i = [1,c_i,d_i]. 
\] 
Then $\underline{B} = (0,1,\infty,b_1,\dots,b_5)$, where the rational functions 
$b_i$ are easy to calculate. The Jacobian determinant 
$|\partial(b_1,\dots,b_5)/\partial(a_1,\dots,a_5)|$ is not identically zero, hence 
it is not zero for general $c_i,d_i$. This implies that the 
image of $\omega(\up)$ must be dense in $Y$. \qed 

\begin{Theorem} 
The variety $\Sigma(5,2,3,8)$ has codimension at least one in $G(3,S_2)$. 
$(N_1 = 55, N_2 = 54).$ \label{theorem.5238} \end{Theorem} 

The machine computation shows that the codimension is exactly one, but we have 
not been able to prove this. 

\demo 
Let $z_0,\dots,z_5$ be the coordinates on $\P^5$. 
Consider the matrix $\left[ \begin{array}{ccc} z_0 & z_1 & z_2 \\ 
z_3 & z_4 & z_5 \end{array} \right]$ and its minors 
\[ G_0 = z_1z_5 - z_2z_4, \; G_1 = z_2z_3 - z_0z_5, 
   \; G_2 = z_0z_4 - z_1z_3. \] 
The locus $G_0 = G_1 = G_2$ is the Segre threefold 
$\sigma(\P^1 \times \P^2)$. 

For $[a,b,c] \in \P^2$, the quadric $a \, G_0 + b \, G_1 + c \, G_2$ is of 
rank $4$, and singular exactly 
along the line joining the points $[a,b,c,0,0,0], [0,0,0,a,b,c]$. 
Denote this line by ${\mathbb M}_{[a,b,c]}$.

Choose general points $p_1, \dots, p_8 \in \P^2$ and 
$Q_1, \dots, Q_8 \in \P^5$. By the lemma, there are points $A_1, \dots, A_8 \in \P^1$
such that $\omega(\up)(\underline{A})$ is the Gale transform of $\uQ$. Hence we may 
as well assume that $Q_i = \sigma(A_i \times p_i)$, i.e., $Q_i \in {\mathbb M}_{p_i}$.  

Let $\Gamma$ be the net 
$\{ [a,b,c] \in \P^2: a \, G_0 + b \, G_1 + c \, G_2 \}$, and define 
\[ \gamma: \P^2 \stackrel{\sim}{\lra} \Gamma, \quad 
 [a,b,c] \lra a \, G_0 + b \, G_1 + c \, G_2. \] 
By construction, $\gamma(p_i)$ is singular at $Q_i$, 
hence $(\Gamma, \emptyset, \gamma)$ is a grove. 
\qed 

\begin{Remark} \rm 
We have failed to produce a geometric argument for the generic uniqueness of 
the grove. But for what it is worth, we can confirm this point 
computationally by reducing the question to linear algebra. 
Choose general points $\up, \uQ$ as above and let 
\[ u_i = a_{i,0} \, z_0^2 + \dots a_{i,20} \, z_5^2, \quad i=0,1,2; \] 
be three quadratic forms, where the $\underline{a}$ are indeterminates. 
Write $p_j = [p_{0,j},p_{1,j},p_{2,j}]$ for $j = 1,\dots,8$. 
Now consider the system 
\[ \begin{aligned}
u_i(Q_j) & = 0 \quad \text{for $j=1,\dots,8$ and $i =0,1,2$.} \\ 
(\sum\limits_{i=0}^2 p_{i,j}\frac{\partial u_i}{\partial z_k}) (Q_j) & = 0 \quad 
\text{for $j=1,\dots,8$ and $k = 0,\dots,5$.}
\end{aligned} \] 
These are $72$ linear homogeneous equations in the $63$ variables $\underline{a}$. 
A Maple calculation shows that for general $\up, \uQ$, there is a 
unique nontrivial solution upto scalars. This solution defines the 
grove for $\up, \uQ$. 
\end{Remark}

\section{Questions} 
In this area, the open problems are certainly not in short supply. 
However, there are four specific themes which we find especially 
appealing. 

1. One would like to have an analogue of the Alexander--Hirschowitz theorem, at least for 
a reasonably broad range of $(n,d,r,s)$. In \cite{Terra1}, 
Terracini claims the following result: 

Assume $n=r=2, d \ge 4$ and $s \ge (d^2+3d+2)/4$ (this is the bound in 
(\ref{bound.dense})). Then $\Sigma$ is dense in $G(2, S_d)$. 

We do not understand his proof and a clarification would be welcome. 

2. Since the imbedding $\Sigma \subseteq G(r, S_d)$ is $GL_{n+1}$ equivariant, 
the equations defining the closure of $\Sigma$ in $G$ are in principle 
expressible in the language of 
classical invariant theory. For small values, there are results making these 
equations explicit. For instance, in the case $(2,3,1,3)$ the hypersurface 
${\overline \Sigma} \subseteq \P^9$ is defined by the Aronhold invariant of ternary 
cubics. Toeplitz \cite{Toeplitz} gives such a combinant for $(3,2,3,5)$, which turns out 
to be a Pfaffian. One would like to have some general theoretical machinery for 
such problems. 

3. Given $\Lambda \in G(r,S_d)$, the locus $\pi_2(\pi_1^{-1}(\Lambda))$ (as 
defined in the introduction) is called the variety of its polar $s$-hedra. It has a 
very rich geometry, see e.g.~\cite{DolgachevKanev,IlRane, RanestadSch} for some old and new 
results. If $n=1$, then it is 
an open subset of a projective space (see \cite{CarSch}), but much remains unknown 
for more than two variables. 

4. We need interesting examples where the class of 
$\overline{\Sigma}$ in the cohomology ring $H^*(G,{\mathbf Z})$ 
can be calculated. For $n=1$, such calculations can be done using the Porteous 
formula (see \cite{ego.Ger}) but in general it is not clear how to proceed. 

\bibliographystyle{plain}
\bibliography{../BIB/hema1,../BIB/hema2,../BIB/hema3}

\bigskip 

\parbox{13cm}{\small 
Enrico Carlini \\
Dipartimento di Matematica \\
Universita' degli Studi Pavia \\ 
Via Ferrata 1 \\ 
27100 Pavia, Italy. \\ 
email: {\tt carlini@dimat.unipv.it}}

\bigskip 

\parbox{13cm}{\small 
Jaydeep V. Chipalkatti \\
Department of Mathematics and Statistics, \\
416 Jeffery Hall, Queen's University, \\ 
Kingston, ON K7L 3N6, Canada. \\ 
email: {\tt jaydeep@mast.queensu.ca}}

\end{document}